\documentclass[12pt,oneside]{article}

\usepackage{amsfonts}
\usepackage{amscd}
\usepackage{amsthm}
\usepackage{amssymb}
\usepackage{amsmath}
\usepackage{epsfig}

\usepackage{graphicx}
\usepackage{xypic}
\input xypic
\input xy
\xyoption{all}

\newcommand{\mD}{\mathbb D}
\newcommand{\mR}{\mathbb R}
\newcommand{\mZ}{\mathbb Z}
\newcommand{\mI}{\mathbb I}
\newcommand{\mH}{\mathbb H}
\newcommand{\mN}{\mathbb N}
\newcommand{\CH}{\check H}
\newcommand{\tH}{\tilde H}

\title{Strong Jordan separation and applications to rigidity.}

\author{Jean-Fran\c cois Lafont}

\theoremstyle{definition}
\newtheorem{Def}{Definition}[section]

\theoremstyle{proposition}
\newtheorem{Lem}{Lemma}[section]
\newtheorem{Prop}{Proposition}[section]

\newtheorem{Claim2}{Claim}

\theoremstyle{plain}

\newtheorem{Thm}{Theorem}[section]
\newtheorem{Cor}{Corollary}[section]

\theoremstyle{remark}

\newtheorem*{Prf}{Proof}

\begin{document}

\maketitle

\begin{abstract}
In this paper, we extend the results of \cite{la} to higher dimension.  We prove
that simple, thick hyperbolic P-manifolds of dimension $\geq 3$ exhibit Mostow
rigidity.  We also prove a quasi-isometry rigidity result for the fundamental
groups of simple, thick hyperbolic P-manifolds of dimension $\geq 3$.  The key tool
in the proof of these rigidity results is a strong form of the Jordan 
separation theorem, for maps from $S^n\rightarrow S^{n+1}$ which are not 
necessarily injective. 
\end{abstract}

\section{Introduction.}

\subsection{Basic definitions.}

In recent years, there has been much interest in proving rigidity type 
theorems for non-positively curved spaces.  All of these results originated
from Mostow's seminal work, in which he showed that homotopy equivalent 
compact locally-symmetric spaces of rank one are always isometric \cite{mo}.

In this paper, we are interested in showing rigidity results for certain
singular CAT(-1) spaces that are ``mostly'' hyperbolic manifolds.  Let us 
start by defining the various objects we will be discussing
in the paper, and state the main theorems we will prove.

\begin{Def}
We define a closed $n$-dimensional {\it piecewise manifold}
(henceforth abbreviated to P-manifold) to be a topological
space which has a natural stratification into pieces which are
manifolds. More precisely, we define a $1$-dimensional P-manifold
to be a finite graph.
An $n$-dimensional P-manifold ($n\geq 2$) is defined inductively as a 
closed pair
$X_{n-1}\subset X_n$ satisfying the following conditions:

\begin{itemize}
\item Each connected component of $X_{n-1}$ is either an
$(n-1)$-dimensional P-manifold, or an $(n-1)$-dimensional manifold.
\item The closure of each connected component of $X_n-X_{n-1}$
is homeomorphic to a compact orientable $n$-manifold with
boundary, and the homeomorphism takes the component of $X_n-X_{n-1}$ 
to the interior of the $n$-manifold with boundary; the closure of such a 
component will be called a {\it chamber}.
\end{itemize}

\noindent Denoting the closures of the connected components of $X_n-X_{n-1}$ 
by $W_i$,
we observe that we have a natural map $\rho: \coprod
\partial (W_i)\longrightarrow X_{n-1}$ from the disjoint union
of the boundary components of the chambers to the subspace
$X_{n-1}$.  We also require this map to be surjective, and a
homeomorphism when restricted to each component. The P-manifold is
said to be \emph{thick} provided that each point in $X_{n-1}$ has
at least three pre-images under $\rho$. We will
henceforth use a superscript $X^n$ to refer to an $n$-dimensional
P-manifold, and will reserve the use of subscripts $X_{n-1},\ldots
,X_1$ to refer to the lower dimensional strata.  For a thick
$n$-dimensional P-manifold, we will call the $X_{n-1}$ strata the
{\it branching locus} of the P-manifold.
\end{Def}

Intuitively, we can think of P-manifolds as being ``built'' by
gluing manifolds with boundary together along lower dimensional
pieces. Examples of P-manifolds include finite graphs and soap 
bubble clusters.  Observe that compact manifolds can also be 
viewed as (non-thick) P-manifolds.
Less trivial examples can be constructed more or less arbitrarily
by finding families of manifolds with homeomorphic boundary and
glueing them together along the boundary using arbitrary
homeomorphisms. We now define the family of metrics we are
interested in.

\begin{Def}
A Riemannian metric on a 1-dimensional P-manifold (finite graph)
is merely a length function on the edge set.  A Riemannian metric
on an $n$-dimensional P-manifold $X^n$ is obtained by first
building a Riemannian metric on the $X_{n-1}$ subspace, then
picking, for each chamber $W_i$ a Riemannian metric with totally
geodesic boundary satisfying that the gluing map $\rho$ is an isometry. 
We say that a Riemannian metric on a P-manifold is
hyperbolic if at each step, the metric
on each $W_i$ is hyperbolic.  
\end{Def}

A hyperbolic P-manifold $X^n$ is automatically a locally $CAT(-1)$ 
space (see Chapter II.11 in Bridson-Haefliger \cite{bh}).  Furthermore, the lower 
dimensional 
strata $X_i$ are totally geodesic subspaces of $X^n$.  In particular,
the universal cover $\tilde X^n$ of a hyperbolic P-manifold $X^n$ is a 
$CAT(-1)$ space (so is automatically $\delta$-hyperbolic), and 
has a well-defined boundary at infinity $\partial
^\infty \tilde X^n$. Finally we note that the fundamental group 
$\pi_1(X_n)$ is a $\delta$-hyperbolic group.  We refer the reader to 
\cite{bh} for background on $CAT(-1)$ and $\delta$-hyperbolic spaces.

\begin{Def}
We say that an $n$-dimensional P-manifold $X^n$ is {\it simple} 
provided its codimension two strata is empty.  In other words, the
$(n-1)$-dimensional strata $X_{n-1}$ consists of a disjoint union
of $(n-1)$-dimensional manifolds.
\end{Def}

Next we introduce a locally defined topological invariant.
We use $\mD^n$ to denote a closed $n$-dimensional disk, and
$\mD^n_\circ$ to denote its interior.  

\begin{Def}
Define the 1-{\it tripod} $T$ to be the topological space obtained
by taking the join of a one point set with a three point set.
Denote by $*$ the point in $T$ corresponding to the one point set.
We define the $n$-{\it tripod} ($n\geq 2$) to be the space $T
\times \mD^{n-1}$, and 
call the subset $*\times \mathbb D^{n-1}$ the {\it spine} of the 
tripod $T\times \mathbb D^{n-1}$.  The subset $*\times 
\mathbb D^{n-1}$ separates
$T\times \mathbb D^{n-1}$ into three open sets, which we 
call the {\it open leaves} of the tripod.  The union of an 
open leaf with the spine will be called a {\it closed leaf}
of the tripod.
We say that a point $p$ in a topological space $X$ is {\it
$n$-branching} provided there is a topological embedding $f:T\times 
\mathbb D^{n-1}
\longrightarrow X$ such that $p\in f(*\times 
\mathbb D^{n-1}_\circ)$.  
\end{Def}

It is clear that the property of being $n$-branching is invariant
under homeomorphisms.  Note that, in a simple, thick P-manifold of 
dimension
$n$, points in the codimension one strata are automatically 
$n$-branching.  One can ask whether this property can be detected
at the level of the boundary at infinity.  This motivated the author \cite{la}
to make the following:

\vskip 5pt

\noindent
{\bf Conjecture:}  Let $X^n$ be a simple, thick hyperbolic
P-manifold of dimension $n$, and let $p$ be a point in the boundary 
at infinity of
$\tilde X^n$.  Then $p$ is $(n-1)$-branching if and only if $p=\gamma
(\infty)$ for some geodesic ray $\gamma$ contained entirely in a
connected lift of $X_{n-1}$.  
\vskip 5pt

Note that in general, the local structure of the boundary at infinity of
a CAT(-1) space (or of a $\delta$-hyperbolic group) is very hard to analyze.  
The conjecture above says that, with respect to branching, the boundary of a
simple, thick hyperbolic P-manifold of dimension $n$ is 
particularly easy to understand.  In \cite{la}, the
author showed that the above conjecture holds for $n=3$.  

\subsection{Statement of results.}

The main goal of this paper is to prove:

\begin{Thm}[Mostow rigidity]
Let $X_1,X_2$ be a pair of simple, thick, n-dimensional ($n\geq 3$), 
hyperbolic P-manifolds,
and assume that $\phi:\pi_1(X_1)\rightarrow \pi_1(X_2)$ is an isomorphism.
Then there is an isometry $f:X_1\rightarrow X_2$ which induces the map 
$\phi$ on the level of fundamental groups.
\end{Thm}   

Note that an immediate consequence of this result is that the canonical 
map $Isom(X)\rightarrow Out(\pi_1(X))$ from the isometry group to the outer
automorphisms of the fundamental group is surjective when $X$ is a simple,
thick, n-dimensional ($n\geq 3$) hyperbolic P-manifold.  Our second result 
determines the coarse geometry of the fundamental groups of these spaces:

\begin{Thm}[Quasi-isometric rigidity]
Let $X$ be a simple, thick, n-dimensional ($n\geq 3$), hyperbolic P-manifold, and
assume that a finitely generated group $\Gamma$ is quasi-isometric to $\pi_1(X)$.  
Then there exists a short exact sequence:
$$0\rightarrow F\rightarrow \Gamma \rightarrow \Lambda \rightarrow 0$$
with $F$ a finite group, and $\Lambda$ a discrete subgroup of $Isom(\tilde X)$.
\end{Thm}

We remark that the 3-dimensional version of these two theorems were proved in the
author's thesis \cite{la}.  In fact, the proofs given in that paper immediately
extend to higher dimension, provided one knows that the conjecture in the
previous section holds.  The author's proof of the 3-dimensional case of
the conjecture in \cite{la}, relied on the fact that the Schoenflies theorem
holds in dimension $2$.  Since the corresponding theorem fails in higher 
dimensions, a different argument is needed.

Now recall that the classical Jordan separation theorem
asserts that the image of an injective map $f:S^n\hookrightarrow S^{n+1}$ 
always separates $S^{n+1}$ into precisely two connected components.
The key result that allows us to extend the 3-dimensional
arguments of \cite{la} to higher dimensions is the following generalization 
of the Jordan separation theorem, which may be of independent interest.

\begin{Thm}[Strong Jordan separation]
Let $f:S^n\rightarrow S^{n+1}$ be a continuous map, and let $I\subset 
S^n$ be the set of injective points (i.e. points $p\in S^n$ with the
property that $f^{-1}(f(p))=\{p\}$).  If $I$ contains an open set $U$,
and $q\in U$, then:
\begin{itemize}
\item $f(S^n)$ separates $S^{n+1}$ into open subsets (we write
$S^{n+1}-f(S^n)$ as a disjoint union $\amalg U_i$, with each $U_i$ open),
\item there are precisely two open subsets $U_1$, $U_2$ in the 
complement of $f(S^n)$ which contain $f(q)$ in their closure.
\item if $F:\mD^{n+1}\rightarrow S^{n+1}$ is an extension of the map $f$ 
to the closed ball, then $F(\mD^{n+1})$ surjects onto either $U_1$ or $U_2$.
\end{itemize}
\end{Thm}

The main application of Strong Jordan separation will be to prove the
conjecture mentioned in the previous section:

\begin{Thm}[Characterization of branching points] 
Let $X^n$ be a simple, thick, $n$-dimensional ($n\geq 3$), hyperbolic
P-manifold, and let $p$ be a point in the boundary 
at infinity of
$\tilde X^n$.  Then $p$ is $(n-1)$-branching if and only if $p=\gamma
(\infty)$ for some geodesic ray $\gamma$ contained entirely in a
connected lift of $X_{n-1}$.  
\end{Thm}

We point out that the quasi-isometry result we have (Theorem 1.2) fits into 
the framework
of understanding quasi-isometries of groups that have a splitting.  In recent
years, there has been an increasing interest in this question (see for instance
\cite{pw}, \cite{kl}, \cite{pap}, \cite{fm1}, \cite{fm2}, \cite{msw}).  Our 
result differs from these in two ways.  Firstly, the groups we allow as
vertex and edge groups differ from those considered in these papers.  Secondly,
our method of proof involves a precise description of the local topology of the 
boundary at infinity (Theorem 1.4).

The outline of the paper is as follows:
in Section 2, we will give a proof of Theorem 1.3 (Strong Jordan separation).
Using Theorem 1.3, we will establish Theorem 1.4 (Characterization of
branching points) in Section 3.  In Section 4, we will give applications of 
Theorem 1.4 by outlining the proof (taken from \cite{la}) of both 
Theorem 1.1 (Mostow rigidity) and Theorem 1.2 
(Quasi-isometric rigidity).  Some concluding remarks will be provided in 
Section 5.

\vskip 10pt 

\centerline{\bf Acknowledgements.}

The author would like to thank his advisor, R. Spatzier, for many helpful 
discussions concerning this project.  He would also like to thank F.T. Farrell
for helpful suggestions on proving Theorem 1.3.  
F.D. Ancel informed the author of the locally-flat 
approximation theorem found in \cite{an},\cite{ac}, \cite{bg} (and used in our proof of
Lemma 2.7).  C. Westerland provided helpful comments on an early draft of this paper.  
Finally the author would like to thank the anonymous referree for pointing out
simplifications for several of the author's original proofs.

\section{Strong Jordan separation theorem.}

In this section, we provide a proof of Theorem 1.3, which will 
subsequently be used in the proofs of the remaining results in this paper.

Before starting with the proof, we note that this theorem clearly 
generalizes the classical Jordan separation theorem (corresponding to
the case $I=S^n$).  The author does not know whether the hypotheses on
$I$ can be weakened to just assuming that $I$ is measurable.  Throughout
our argument, $\CH^*$ will refer to \v{C}ech cohomology with coefficients in 
$\mZ$, while $\tilde H_*$ will refer to reduced singular homology with 
coefficients in 
$\mZ$.  We use $\mD^n$ to refer to a closed $n$-dimensional disk, and 
$\mD^n_\circ$ to refer to its interior.  We also occasionaly use $\mI$
to denote a closed interval, and $\mI_\circ$ for its interior.

\subsection{Global separation.}

In this subsection, we provide a proof of the first claim in Theorem 1.3, and
some other results which are of a global nature.

\begin{Lem}
Under the hypotheses above, $f(S^n)$ separates $S^{n+1}$.
\end{Lem}

\begin{Prf}
We start by showing that $f(S^n)$ cannot surject onto $S^{n+1}$.  Since $I$
is assumed to contain an open set, we can find a $\mD^n\subset I\subset S^n$,
i.e. a (small) disk on which $f$ is injective.  Let $p\in \mD^n_\circ$ 
be a point in the interior of the disk, and note that (since $p\in I$), we 
have that $d(f(p),f(q))\neq 0$ for all $q\in S^n-\{p\}$.  Since both 
$p$ and $f(S^n-\mD^n_\circ)$ are compact sets, the distance between them must
be positive.  Pick $\epsilon$ smaller than $d(f(p),f(S^n-\mD^n_\circ))$, and note
that an $\epsilon$-ball centered at $f(p)$ can only intersect $f(S^n)$ in a 
subset of $f(\mD^n)$.  But on this set, $f$ is injective, and it is well known
that there is no injective and surjective map from a subset of $\mD^n$ to $\mD^{n+1}$ 
(the $\epsilon$-ball centered at $p$).  This shows that $f(S^n)$ is not surjective.

Since $f$ is not surjective, Alexander duality tells us that 
$\tilde H_0(S^{n+1}-f(S^n))\cong \CH^n(f(S^n))$, hence it is sufficient to 
show that $\CH^n(f(S^n)) \neq 0$.   Now the space $f(S^n)$ 
can be expressed as a union $f(S^n)=f(\mD^n)\cup_g f(S^n-\mD^n_\circ)$
(where $g:\partial \mD^n=S^{n-1} \rightarrow f(S^n-\mD^n_\circ)$ is the 
gluing map); this follows from the fact that the defining property of the
set $I$ ensures that $f(\mD^n)$ intersects $f(S^n-\mD^n_\circ)$ in 
precisely $f(\partial \mD^n)$. 
Since $f$ is injective on $\mD^n$, we have $f(\mD^n)$ is homeomorphic to 
$\mD^n$, and the gluing is along $\partial \mD^n = S^{n-1}$.  Furthermore,
note that $g$ is homotopic to
a point in $f(S^n-\mD^n_\circ)$, and since $S^{n-1}$ is collarable in $\mD^n$,
we have that the space $f(S^n)$ is homotopy
equivalent to the space $X=f(\mD^n)\cup_h f(S^n-\mD^n_\circ)$, where $h$ maps 
$\partial \mD^n$ to a point.  In particular, we see that $f(S^n)$ is
homotopy equivalent to the join $S^n\vee f(S^n-\mD^n_\circ)$.  This implies
that $\CH^i(f(S^n))\cong \CH^i(S^n \vee f(S^n-\mD^n_\circ))$.  But we can
compute the latter by using the Mayer-Vietoris sequence in \v{C}ech cohomology,
yielding $\CH^n(f(S^n))\cong \mZ \oplus \CH^n(f(S^n-\mD^n_\circ))\neq 0$.  This
gives us our claim.
\end{Prf}

We now have that the image does indeed separate $S^{n+1}$.  We denote by
$U_i$ the connected components of $S^{n+1}-f(S^n)$, giving us the decomposition
$S^{n+1}-f(S^n)=\amalg U_i$.  Note that $S^{n+1}-f(S^n)$ is open, so that
connected components coincide with path-connected components.  Let us now
focus on the local behavior near $f(\mD^n_\circ)$.

\begin{Lem}
The inclusion map $S^{n+1}-f(S^n)
\hookrightarrow S^{n+1}-f(S^n-\mD^n_\circ)$ induces a canonical splitting
$\tilde H_0(S^{n+1}-f(S^n))\cong \mZ \oplus 
\tilde H_0(S^{n+1}-f(S^n-\mD^n_\circ))$.
\end{Lem}

\begin{Prf}
Since Alexander duality is natural, we have a commutative diagram:
$$\xymatrix{\tilde H_0(S^{n+1}-f(S^n)) \ar[r] \ar[d]^\cong & 
\tilde H_0(S^{n+1}-f(S^n-\mD^n_\circ))
\ar[d]^\cong \\ \CH^n(f(S^n)) \ar[r]& \CH^n(f(S^n-\mD^n_\circ))}
$$
where the two horizontal arrows are induced by the obvious inclusions,
and the vertical arrows are given by duality.  So it is sufficient to
show that the morphism $\CH^n(f(S^n)) \rightarrow \CH^n(f(S^n-
\mD^n_\circ))$ gives the desired splitting. But this is clear, since 
we have $\CH^n(f(S^n))\cong \mZ \oplus \CH^n(f(S^n-\mD^n_\circ))$, where the 
direct sum comes from the Mayer-Vietoris sequence, and
the $\mZ$ comes from $\CH^n(S^n)$ (see previous Lemma).  Combining this 
with the commutative diagram above gives us the desired splitting
$\tilde H_0(S^{n+1}-f(S^n))\cong \mZ \oplus 
\tilde H_0(S^{n+1}-f(S^n-\mD^n_\circ))$, completing 
the proof of the lemma.
\end{Prf}

This gives us a homological version of the statement ``$f(\mD^n_\circ)$ 
locally seperates $S^{n+1}$ into pieces lying in precisely two components
of $S^{n+1}-f(S^n)$''.  We now want to
replace this by a non-homological, equivalent statement.

\begin{Lem}
Let $\{U_i\}$ be a collection of disjoint connected open sets in $\mR^{n+1}$.
Assume that there is a connected set $Z$ which intersects the closure of 
each $U_i$ non-trivially.  Then the space $Y:=Z\cup (\bigcup U_i)$ is a 
connected set. 
\end{Lem} 

The proof of Lemma 2.3 is an easy exercise in point-set topology, which we 
leave to the reader.

%

\begin{Cor}
There are precisely two connected components $U_1$, $U_2$ in $S^{n+1}-
f(S^n)$ whose closure intersects $f(\mD^n_\circ)$. 
\end{Cor}

\begin{Prf}
Let $r$ denote the number of connected components in $S^{n+1}-f(S^n)$
whose closure intersects $f(\mD^n_\circ)$, and let us enumerate them
by $\{U_i\}_{i=1}^r$.  Picking points $x_i\in U_i$, note that for 
$1\leq j<k\leq r$, all the elements $[x_j-x_k]\in \tilde H_0(S^{n+1}-f(S^n))$
represent distinct non-zero classes.  On the other hand, Lemma 2.3 tells us
that each $U_j,U_k$ must lie in the same connected component of 
$S^{n+1}-f(S^n-\mD^n_\circ)$, and hence that the classes $[x_j-x_k]$
map to zero under the homomorphism:
$$\tilde H_0(S^{n+1}-f(S^n))\rightarrow \tilde H_0(S^{n+1}-f(S^n-\mD^n_\circ))$$
induced by the inclusion.  This implies that the kernel of the 
homomorphism is in fact generated by the classes $[x_i-x_{i+1}]$ ($1\leq i<r$).

But now Lemma 2.2 tells us that the homomorphism above induces a 
splitting $\mZ \oplus \tilde H_0(S^{n+1}-f(S^n-\mD^n_\circ)) \cong \tilde 
H_0(S^{n+1}-f(S^n))$.  Since the kernel has rank 1, this implies that
there are precisely two connected components in $S^{n+1}-f(S^n)$ 
whose closure intersects $f(\mD^n_\circ)$, as desired.
\end{Prf}

\subsection{Local separation.}

While the separation results we've obtained so far are useful, they do
not tell us much about the behavior of the embedding near a point in the 
injectivity set $I$.  The purpose of this section is to fine-tune the 
information we have in order to capture local information concerning the
map $f$.  We start with a local version of Lemma 2.2:

\begin{Lem}
Let $p\in f(\mD^n_\circ)$ be an arbitrary point, and $U_1, U_2$ be the two 
connected components of $S^{n+1}-f(S^n)$ whose closure intersects 
$f(\mD^n_\circ)$.  Then $p$ lies in the closure of both $U_1$ and 
$U_2$.
\end{Lem}

\begin{Prf}
Given $p\in \mD^n_\circ$, choose a sequence of open metric balls 
$D_i\subset \mD^n_\circ$
centered at $p$, having the property that their radii tend to zero.  
Now since each $D_i$ lies
in the injectivity set $I$, Corollary 2.1 tells us that for each $D_i$,
there are precisely two connected components $U_1^{(i)}, U_2^{(i)}$
in $S^{n+1}-f(S^n)$ whose closure intersects $f(D_i)$.  But since 
$f(D_i)\subset f(\mD^n_\circ)$, this forces (upto relabelling) 
$U_1^{(i)}=U_1$, $U_2^{(i)}=U_2$.  This implies that the closure
of $U_1,U_2$ intersects points in each $f(D_i)$.  Since the radii of the
$D_i$ tends to zero, this forces $f(p)$ to likewise lie in the closure
of both $U_1, U_2$.
\end{Prf}

Next we introduce a group that measures, for an arbitrary point $p$ in 
$S^{n+1}$, the number of components of $S^{n+1}-f(S^n)$ that contain 
$p$ in their closure.

\begin{Def}
For a point $p\in S^{n+1}$, and a continuous map $f:S^n\rightarrow 
S^{n+1}$, we define the {\it group of local components of $p$ relative 
to $f$} to be:
$$H(p,f):= \underleftarrow{\lim} \{\tH _0(B(p,1/m) -f(S^n))\},$$
where $B(p,1/m)$ is the open metric ball of radius $1/m$, centered at $p$, 
and the inverse limit is over the directed system $m\in \mN$, with 
morphisms induced by inclusions.
\end{Def}

Observe that the groups $H(p,f)$ are {\it local}, in the sense that if 
one has a pair of maps $f,g:S^n\rightarrow S^{n+1}$, and the point $p$ 
has an open neighborhood $N$ with the property that $f(S^n)\cap 
N=g(S^n)\cap N$, then the groups $H(p,f)$ and $H(p,g)$ are canonically 
isomorphic.  Furthermore, note that if $p \not \in f(S^n)$, then 
$H(p,f)= 0$.

Next we observe that for an arbitrary point $p\in S^{n+1}$, and each 
$m\in \mN$, we have inclusions $B(p,1/m)-f(S^n)\hookrightarrow 
S^{n+1}-f(S^n)$, which induce homomorphisms $\tH _0(B(p,1/m) -f(S^n)) 
\rightarrow \tH_0(S^{n+1}-f(S^n))$.  It is easy to see that these 
homomorphisms induce a homomorphism $\phi:H(p,f)\rightarrow 
\tH_0(S^{n+1}-f(S^n))$.  Our next two lemmas serve to analyze the image 
of this map in the case where $p\in f(\mD^n_\circ)$.

\begin{Lem}
If $p\in f(\mD^n_\circ)$, the $H(p,f)\cong \mZ$.
\end{Lem}

\begin{Prf}
Fix $m\in \mN$, and consider the
map $H:S^{n+1}\rightarrow S^{n+1}=B(p,1/m)^*$ obtained by collapsing the 
complement of $B(p,1/m)$ to a point (so that the target can be identified with 
the one point compactification $B(p,1/m)^*$ of the open ball $B(p,1/m)$).  
Applying Corollary 2.1 and Lemma 2.4 to the map 
$H\circ f:S^n\rightarrow B(p,1/m)^*$, we see that there are exactly two connected
components of $B(p,1/m)^*-(H\circ f)(S^n)$ containing $p$ in their closure.  
But this immediately implies that there are precisely two connected components 
of $B(p,1/m)-f(S^n)$ that contain $p$ in their closure

Now consider what this tells us about the inverse system of groups 
$\{\tH_0(B(p,1/m)-f(S^n))\}_{m\in \mN}$.  By an argument similar to the one given at 
the beginning of the proof of Lemma 2.4, we see that for each $m\in \mN$,
there exists a $M$ with the property that $\forall k\geq M$, the set 
$B(p,1/k)-f(S^n)$ lies in the union of the two connected components of $B(p,1/m)-
f(S^n)$ which contain $p$ in their closure.  In particular, on the level
of $\tH_0$, we see that for all $k\geq M$, the map $\tH_0(B(p,1/k)-f(S^n))
\rightarrow \tH_0(B(p,1/m)-f(S^n))$ has image lying in a fixed $\mZ$-subgroup
of $\tH_0(B(p,1/m)-f(S^n))$.  But the inverse limit of an inverse system
of groups with this
property is automatically isomorphic to $\mZ$, concluding the proof of
Lemma 2.5.
\end{Prf}

\begin{Lem}
If $p\in f(\mD^n_\circ)$, then $\phi(H(p,f)) \cong \mZ$.  Furthermore, 
the image of $\phi$ is in fact generated by the reduced homology class 
$[x_1-x_2]\in \tH_0(S^{n+1}-f(S^n))$, where $x_i\in U_i$ (and $U_1,U_2$ are
the two connected components of $S^{n+1}-f(S^n)$ whose closures intersect
$f(\mD^n_\circ)$).
\end{Lem}

\begin{Prf}
If $m\in \mN$ is arbitrary, then $B(p,1/m)-f(S^n)$ must contain points in 
both $U_1$ and $U_2$.  Hence the inclusion induces a map $\tH_0(B(p,1/m)-f(S^n))
\rightarrow \tH_0(S^{n+1}-f(S^n))$ whose image contains the subgroup generated
by the reduced homology class $[x_1-x_2]\in \tH_0(S^{n+1}-f(S^n))$.  Since
this holds for all $m\in \mN$, the induced map on the inverse limit 
$H(p,f)\rightarrow \tH_0(S^{n+1}-f(S^n))$ likewise has image containing the
subgroup generated by $[x_1-x_2]$.

Furthermore, if $m$ is sufficiently large, then the image of the map 
$\tH_0(B(p,1/m)-f(S^n))\rightarrow \tH_0(S^{n+1}-f(S^n))$ actually lies in
the subgroup generated by $[x_1-x_2]$.  Passing to the induced map from the
inverse limit $H(f,p)$, we see that $\phi$ maps $H(f,p)$ into the subgroup
generated by $[x_1-x_2]$, completing our proof.
\end{Prf}

Morally speaking, the last two lemmas make precise the fact that 
{\it the components which are incident to a point in the set of 
injectivity can be detected purely locally.}  Indeed, the combination
of Lemma 2.5 and 2.6 tells us that the map $\phi$ is in fact an isomorphism
from $H(p,f)$ to the subgroup generated by $[x_1-x_2]$.  This will be crucial 
to us in the next section, where we prove the third claim from our theorem.

\subsection{Surjectivity of extensions.}

Finally, to complete the proof of the theorem, we need to show that any
extension $F:\mD^{n+1}\rightarrow S^{n+1}$ of the map $f$ must surject onto
either $U_1$ or $U_2$.  Let us assume, without loss of generality, that $F$ does 
not surject onto $U_1$.  Then using a stereographic projection, 
we can view the image of $f$ as lying in $\mR^{n+1}$, with $U_1$ unbounded and all 
other components in the complement of $f(S^n)$ being bounded.  Furthermore, we also
have that the image of $F$ lies in $\mR^{n+1}$ (after stereographic projection).

We will also fix an orientation
on $\mR^{n+1}$, yielding for every point $z\in \mR^{n+1}$ a canonical generator 
for $H_n(\mR^{n+1}-z)\cong \mZ$.  In particular, using this identification,
it makes sense to say that an $n$-cycle in $\mR^{n+1}-\{z_1,z_2\}$ represents 
distinct elements in $H_n(\mR^{n+1}-z_1)$ and $H_n(\mR^{n+1}-z_2)$.

Now observe that for any pair of points $z_1,z_2$ which lie in a fixed connected
component $U$ of $\mR^{n+1}-f(S^n)$, we have that the map $f$ represents the 
same element in $H_n(\mR^{n+1}-z_1)$ and $H_n(\mR^{n+1}-z_2)$.  Furthermore, 
if there exists an extension of the map $f$ to a map $F:\mD^{n+1}\rightarrow \mR^{n+1}$
which fails to surject onto the component $U$, then denoting by $z\in U$ a point
that does not lie in the image of $F$, we see that $f$ represents $0$ in the group
$H_n(\mR^{n+1}-z)$ (since $f$ bounds in the complement of $z$), and hence that $f$
represents $0$ in every $H_n(\mR^{n+1}-z^\prime)$ where $z^\prime \in U$ is arbitrary.
In particular, since the image of $F$ is compact, and $U_1$ is unbounded, we have
that $f$ must represent $0$ in $H_n(\mR^{n+1}-z^\prime)$ whenever $z^\prime \in U_1$.

With this in mind, the last claim of the theorem will follow immediately from:

\begin{Lem}
The map $f:S^n\rightarrow \mR^{n+1}$, ($n\geq 2$) cannot represent $0$ in 
$H_n(\mR^{n+1}-z)\cong \mZ$, where $z\in U_2$. In particular, $F$ must surject onto $U_2$.
\end{Lem}

\begin{Prf}
Note that we may assume that the set $\mD^n\subset
S^n$ is actually a round disk; let $p\in I$ be the center of $\mD^n$, $r$
the radius of $\mD^n$.  Denote by  $\hat \mD^n$ the closed round disk with 
center $p$, and radius $r/2$, and for $t\in [0,r/2]$, let $d(t)$ be the 
distance from $f(S^{n-1}_t)$ to $f(S^n-\hat D^n_\circ)$, where $S^{n-1}_t$ is 
the sphere of radius $t$ centered at $p$.  Note that, for each $t\in [0,r/2]$, 
the sphere
$S^{n-1}_t$ lies in $\mD^n$, hence is contained in the set of injective points.
Let $q\in \mD^n_\circ-\hat \mD^n$ be an arbitrary point, and let
$[x_m-y_m] \in \tH_0(B(f(q),1/m)-f(S^n))$ be a sequence of group elements 
defining the generator in $H(f(q),f)$ which maps to $[z^\prime -z]\in \tH_0(S^{n+1}
-f(S^n))$ (recall that $z^\prime \in U_1$ and $z\in U_2$).   

We would now like to approximate the map $f$ by a new map $g$ which has the
following properties:
\begin{itemize}
\item $f=g$ on the complement of $\hat \mD^n_\circ$,
\item $g$ is locally flat embedding on $\hat \mD^n_\circ$,
\item $g$ is an embedding on $\mD^n_\circ$
\end{itemize}
To see that such an approximation exists, we take a continuous function 
$\delta:[0,r/2)\rightarrow \mR^+$ with the property that $\delta (t)< d(t)$.  
A result of Ancel-Cannon \cite{ac} ($n\geq 4$), 
Ancel \cite{an} ($n=3$), and Bing \cite{bg} ($n=2$), states that
we can approximate the map $f$ on $\hat \mD^n_\circ$, by a map $\hat g:\hat \mD^n_\circ
\rightarrow \mR^{n+1}$ having the property that $\hat g$
is a locally flat embedding, and $d(\hat g, f)\leq \delta(t)$ in the sup norm on
each $S^{n-1}_t$ ($t\in [0, r/2)$).  Finally, we note that the concatenation of 
$\hat g$ on the $S^{n-1}_t$ ($t\in [0,r/2)$) with the map $f$ on the $S^{n-1}_t$ 
($t\in [r/2,r]$) has the desired properties; let us call this map $g$.

Applying Corollary 2.1 to the map $g$, let us  
denote by $U_1^\prime$ and $U_2^\prime$ the two components of $S^{n+1}-g(S^n)$
that contain $g(p)$ in their closure.  From Lemma 2.4, we know that $f(q)=g(q)$ also
lies in the closure of both $U_1^\prime$ and $U_2^\prime$.
Next we note that, for a small enough neighborhood $N$ of $f(q)$, we have that 
$N\cap f(S^n)=N\cap g(S^n)$, hence $H(f(q),f)$ is canonically isomorphic to 
$H(f(q),g)$.  In particular, the sequence $[x_m-y_m]\in \tH_0(B(f(q),1/m)-f(S^n))$
still defines a generator for $H(f(q),g)$.  This implies that, for $m$ large 
enough $x_m\in U_1^\prime$ and $y_m\in U_2^\prime$ (upto reindexing).

Now observe that, since $g$ is locally flat near $\hat \mD^n_\circ$, the map
$g$ {\it is actually flat} on $\hat \mD^n$.  A proof of this old folklore result
can be found in Rushing \cite{ru} Flattening Theorem 3.4.1.
This implies that there
is a homeomorphism $\phi:\mR^{n+1}\rightarrow \mR^{n+1}$ with the property
that $(\phi\circ g)(\hat \mD^n_\circ)\subset \mR^n\times \{0\}$.  Furthermore,
this homeomorphism takes the $U_i^\prime$ to the two open sets $W^+$ 
(respectively $W^-$) containing $(\phi\circ g)(\hat \mD^n_\circ)$ in their
closure, where $W^\pm$ consists of points lying (locally near 
$(\phi\circ g)(\hat \mD^n_\circ)$) above/below the hyperplane 
$\mR^n\times \{0\}$.  But now it is classical that for a pair of points 
$x\in U_1^\prime$, $y\in U_2^\prime$, that the map $\phi\circ g$ must represent 
distinct elements
in $H_n(\mR^{n+1}-\phi(x))$ and $H_n(\mR^{n+1}-\phi(y))$.  Indeed, there is 
a PL-curve
joining $\phi(x)$ to $\phi(y)$ and intersecting $\phi\circ g$ transversally 
in a single point.  Pulling back along $\phi$, we see that $g$ represents
distinct elements in $H_n(\mR^{n+1}-x)$ and $H_n(\mR^{n+1}-y)$.  

This tells us that for $m$ large enough, the map $g$ represents distinct
elements in $H_n(\mR^{n+1}-x_m)$ and $H_n(\mR^{n+1}-y_m)$.  Finally, note
that by our choice of the control function $\delta$, for $m$ large enough
we have that $g$ is homotopic to $f$ in the complement of $\{x_m,y_m\}$.
Hence $f$ must represent distinct elements in $H_n(\mR^{n+1}-x_m)$ and 
$H_n(\mR^{n+1}-y_m)$.  Since $f$ represents $0$ in $H_n(\mR^{n+1}-x_m)$
(by the comment before this proof), we get the desired result. 
\end{Prf}

Combining Lemma 2.1, Corollary 2.1, and Lemma 2.7 proves Theorem 1.3.

\section{Recognizing branching.}

In this section, we will provide a proof of Theorem 1.4 from the
introduction.  For the convenience of the reader, we restate the
Theorem:

\vskip 5pt

\noindent 
{\bf Theorem 1.4.}
{\it Let $X^n$ be a simple, thick hyperbolic
P-manifold of dimension $n\geq 3$, and let $p$ be a point in the boundary 
at infinity of
$\tilde X^n$.  Then $p$ is $(n-1)$-branching if and only if $p=\gamma
(\infty)$ for some geodesic ray $\gamma$ contained entirely in a
connected lift of $X_{n-1}$ (the $(n-1)$-dimensional strata).}

\vskip 5pt

As mentioned in the introduction, the 3-dimensional version of this result
was shown by the author in \cite{la}.  The key point where the hypothesis $n=3$
was used was in the 2-dimensional version of the Strong Jordan Separation
Theorem (Theorem 2.1 in \cite{la}).  Indeed, in that proof, we made use of the
Schoenflies Theorem, which holds in dimension $2$, but is known to fail in
dimension $\geq 3$.  Since the majority of the results in this section 
follow immediately from the proofs of the corresponding results in \cite{la},
we will merely emphasize the arguments that need to be modified.  We start 
by mentioning:

\begin{Prop}
Let $X^n$ be a simple, thick $n$-dimensional 
P-manifold , and let $p$ be a point in the boundary at infinity of
$\tilde X^n$. If $\gamma$ is a geodesic ray contained entirely in a
connected lift $\tilde B$ of $X_{n-1}$, then $\gamma(\infty)$ is 
$(n-1)$-branching.
\end{Prop}

This was Proposition 2.1 in \cite{la}.  The idea was to use the simplicity 
and thickness
hypothesis to ensure that there are three disjoint totally geodesic subsets,
each isometric to a ``half'' $\mH^n$, whose boundaries all coincide with
$\tilde B$.  This forces all points in $\partial ^\infty \tilde B$ to be
$(n-1)$-branching.

Note that points in $\partial ^\infty \tilde X^n$ which {\it do not} lie
in any of the $\partial ^\infty \tilde B_i$ are of two types.  They are
either:
\begin{itemize}
\item endpoints of geodesic rays that eventually stay trapped in a connected
lift $\tilde W$ of a chamber $W$, and are {\it not} asymptotic to any boundary
component of $\tilde W$, or
\item endpoints of geodesic rays that pass through infinitely many connected
lifts $\tilde W_i$ of chambers $W_i$.
\end{itemize}
We will discuss these two cases separately.

\subsection{The geodesic ray is trapped in a chamber.}

We now explain how to deal with the first of the two cases.

\begin{Prop}
Let $X^n$ be a simple, thick $n$-dimensional hyperbolic P-manifold ($n\geq 3$).
Let $\gamma\subset \tilde X^n$ be a geodesic ray lying entirely
within a connected lift $\tilde W$ of a chamber $W$, and not
asymptotic to any boundary component of $\tilde W$.  Then
$\gamma (\infty)$ is {\bf not} $(n-1)$-branching.
\end{Prop}

This was Proposition 2.2 in \cite{la}.  The argument relies on picking a point
$x\in \gamma$, and considering the link of the point $x$ in $\tilde X^n$.
Note that since $x$ is in a chamber, the link $lk(x)$ is isometric to
$S^{n-1}$.  We denote by $\pi_x:\partial^\infty \tilde X^n\rightarrow lk(x)
\cong S^{n-1}$ the geodesic retraction from the boundary at infinity to 
the link at $x$.  Finally, we denote by $I_x$ the set $\{p\in
lk(x) \hskip 5pt : \hskip 5pt |\pi_x^{-1}(p)|=1\}\subset
lk(x)$, in other words, the set of points in the link where
the projection map is actually injective.

Now by way of contradiction, assume that $\gamma (\infty)$ is $(n-1)$-branching.
Then there exists an injective map $f:T\times \mD^{n-2}\rightarrow \partial
^\infty \tilde X^n$ with the property that $\gamma (\infty)\in f(*\times
\mD^{n-1}_\circ)$.  One can now consider the composite map $\pi_x \circ f:
T\times \mD^{n-1}\rightarrow lk(x)$.  Since $lk(x)$ is topologically $S^{n-1}$,
we know that the composite map cannot be an embedding.  The strategy lies
in showing that there exists a point in $I_x$ which has two distinct pre-images
under the composite $\pi_x\circ f$.  Since $\pi_x$ is injective on points in
$I_x$, this would imply that $f$ fails to be injective, yielding a contradiction.

To see that the composite $\pi_x\circ f$ fails to be injective at a point in 
$I_x$, we proceed via a series of Claims.

\begin{Claim2}
The complement of $I_x$ has the following properties:
\begin{itemize}
\item it consists of a countable union of {\it open} disks 
$U_i$ in $S^{n-1}$,
\item the $U_i$ are the interiors of a family of pairwise 
disjoint closed disks,
\item the $U_i$ are dense in $S^{n-1}$.
\end{itemize}
\end{Claim2}

\begin{Claim2}
For any point
$p\in I_x -\cup (\partial U_i)$, and any neighborhood $N_p$ of $p$, 
there exist arbitrarily small $U_i$ with $U_i\subset N_p$.
\end{Claim2}

\begin{Claim2}
The image $(\pi_x\circ f)(\partial(T\times \mD^{n-2}))$ is a bounded
distance away from $\pi_x(\gamma(\infty))$.
\end{Claim2}

The proof of these three claims given in \cite{la} for the case $n=3$ extend 
verbatim to the case $n\geq 3$.  For the convenience of the reader, we briefly
summarize the arguments.

Points in the complement of $I_x$ correspond to directions in which $\pi_x$ 
fails to be injective.  This means that there exists a pair of geodesic rays
which coincide on a small neighborhood of $x$, and subsequently diverge.  So
the complement of $I_x$ is precisely the image of $\partial \tilde W$ under the 
geodesic retraction.  Since $\tilde W$ is the universal cover of a compact
hyperbolic manifold with boundary, standard methods yield the first claim.

For the second claim, we use the fact that the $U_i$ are 
actually {\it round} disks in $lk(x)=S^{n-1}$.  This implies that, for the 
standard measure $\mu$ on $S^{n-1}$, the radius of the $U_i$ is proportional to
$\mu(U_i)^{\frac{1}{n-1}}$.  In particular, there are only finitely many $U_i$
of radius $\geq r$ for any $r>0$.  Using this fact, one can
find a suitably smaller neighborhood $N_p^\prime\subset N_p$ with the property
that the $U_i$ that intersect $N_p^\prime$ have radius less than half the 
distance between $\partial N^\prime_p$ and $S^{n+1}- N_p$, giving the
second claim.

The third claim follows immediately from the compactness of the set 
$(\pi_x\circ f)(\partial(T\times \mD^{n-2}))$, along with the fact that 
$\pi_x(\gamma(\infty))\notin (\pi_x\circ f)(\partial(T\times \mD^{n-2}))$.
\vskip 10pt

Continuing with our argument we now have:

\begin{Claim2}
There is a neighborhood $N$ of the point $f^{-1}(\gamma (\infty))$ in
the spine $*\times \mD^{n-2}_\circ$ whose image is entirely contained in $I_x$ 
(i.e. $(\pi_x\circ f)(N)\subset I_x$).
\end{Claim2}

This claim requires a substantial modification from that given in \cite{la}.
We first note that, denoting by $r_i$ the radius of the respective $U_i$, 
the argument for Claim 2 shows that for any $R>0$, there are only finitely 
many $U_i$ with the property that $r_i>R$.  Note that, by hypothesis, we
have that $\pi_x(\gamma (\infty))\notin \cup \bar U_i$.  Hence it is 
sufficient to show that, for a small enough $\epsilon$, there are only
finitely many $U_i$ intersecting $(\pi_x\circ f)(*\times \mD^{n-2}_\circ)
\cap B(\pi_x(\gamma(\infty)))$, where $B(\pi_x(\gamma(\infty)))$ denotes
the open $\epsilon$ ball centered at $\pi_x(\gamma(\infty))$.

Now if this were not the case, then there exists, for each $\epsilon >0$,
infinitely many $U_i$ intersecting $(\pi_x\circ f)(*\times \mD^{n-2}_\circ)$
at a distance $<\epsilon$ from $\pi_x(\gamma (\infty))$.  Now pick $\epsilon$
to be strictly smaller than the distance $D$ from $\pi_x(\gamma (\infty))$ 
to the set $(\pi_x\circ f)(\partial(T\times \mD^{n-2}))$.  We can now find 
a $U_i$ such that $2r_i+\epsilon <D$, {\it and} $U_i$ intersects 
$(\pi_x\circ f)(*\times \mD^{n-2}_\circ)$ at a distance $<\epsilon$ from
$\pi_x(\gamma (\infty))$.  

Let us focus now on this $U_i$.  The condition that $2r_i+\epsilon <D$ 
ensures that $\bar U_i \cap (\pi_x\circ f)(\partial(T\times \mD^{n-2}))=
\emptyset$ (i.e. that the image of the boundary of the $(n-1)$-tripod fails to 
intersect the closure of the $U_i$). 
Furthermore, note that $\partial U_i\subset I_x$, and is homeomorphic to 
$S^{n-2}$.  This implies that $(\pi_x\circ f)^{-1}(\partial U_i)$ is a closed
subset of $S^{n-2}$ 
(topologically) embedded in $T\times \mD^{n-2}_\circ$, and which {\it
fails} to intersect $\partial (T\times \mD^{n-2})$.  Since
$\partial U_i$ separates $lk(x)$, we have that $(\pi_x\circ f)^{-1}(\partial
U_i)$ likewise separates $T\times \mD^{n-2}_\circ$.

Let us denote by $S$ the subset $(\pi_x\circ f)^{-1}(\partial U_i)$ in
$T\times
\mD^{n-2}_\circ$, and denote by $L_1, L_2, L_3$ the three closed leaves of 
$T\times\mD^{n-2}_\circ$.  We now observe that $L_j\cup L_k$ is homeomorphic 
to $\mD^{n-1}$ ($1\leq j<k\leq3$), and that the intersection 
$S_{j,k}:=S\cap (L_j\cup L_k)$ is a closed subset of $L_j\cup L_k$. 
Furthermore, we note that at least two of the three sets $S_{1,2}, S_{1,3}, 
S_{2,3}$ must separate the corresponding $L_j\cup L_k$.  Without loss of
generality, assume that $S_{1,2}, S_{2,3}$ separate the corresponding
$L_j\cup L_k$.  Finally, observe that since $S \cap \partial (T\times
\mD^{n-2})=\emptyset$, we can view each of the sets $S_{1,2}$, $S_{2,3}$
as lying in the interior of the corresponding $L_i\cup L_j$, which is
naturally homeomorphic to $\mR^{n-1}$.

Now note that both $S_{1,2}$ and $S_{2,3}$ are closed subsets 
which seperate $Int(L_1\cup L_2)\cong \mR^{n-1}$ and $Int(L_2\cup 
L_3)\cong \mR^{n-1}$ respectively.  By Alexander duality, this implies that
$\CH^{n-2}(S_{1,2})\neq 0$ and $\CH^{n-2}(S_{2,3})\neq 0$.  However, both these
sets are homeomorphic to closed subsets of $S^{n-2}$, and the only closed
subset $A$ of $S^{n-2}$ satisfying $\CH^{n-2}(A)\neq 0$ is $S^{n-2}$ itself.
This implies that both $S_{1,2}$ and $S_{2,3}$ must be homeomorphic to 
$S^{n-2}$, yielding that
$S_{1,2}=S_{2,3}=S$.  From the definition of the sets $S_{j,k}$, we immediately
get that $S\subset L_2$.  

So we have now reduced the problem to showing that the pre-image 
$S=(\pi_x\circ f)^{-1}(\partial U_i)$ of 
$\partial U_i$ cannot lie entirely in $L_2$.  Note that $S$ separates $L_2$
into two components, one of which must map into $U_i$, while the other one
must map into $lk(x)-U_i$.  Now by hypothesis, there is a point in the spine
$*\times \mD^{n-2}$ of $T\times \mD^{n-2}$ that maps into $U_i$, which implies
that there is a point on the boundary of $L_2=\mI\times \mD^{n-2}$ which maps 
into $U_i$.  This implies that the entire spine must in fact map either to
$U_i$, or to $\partial U_i$.  Since $\pi_x(\gamma(\infty))$ lies in the
image of the spine, we contradict the hypothesis that
$\pi_x(\gamma(\infty))\in I_x$ does not lie in the closure of any $U_i$.
This concludes the argument for Claim 4.

\vskip 10pt

Using Claim 4, and focusing on a sufficiently small neighborhood of 
$f^{-1}(\gamma (\infty))$ in $T\times \mD^{n-2}$, we may assume that the
map $\pi_x\circ f$ maps the spine into the set $I_x$.  We will henceforth
work under this assumption.

Let us denote by $F_i$ the restriction of the composition $\pi_x\circ f$
to each of the closed leaves in $T\times \mD^{n-2}$, and by $G_i$ the restriction
of $\pi_x\circ f$ to the boundary of the closed leaves.  We now have the following:

\begin{Claim2}
There is a connected open set $V\subset S^{n-1}$ with the property
that:
\begin{itemize}
\item at least two of the maps $F_i$ surject onto $V$
\item the closure of $V$ contains the point $\pi_x(\gamma (\infty))$
\end{itemize}
\end{Claim2}

As in \cite{la}, Claim 5 is a consequence of the Strong 
Jordan separation theorem.  
Observe that each leaf is homeomorphic to $\mI\times \mD^{n-2} \cong \mD^{n-1}$, 
and from Claim 4, we know that the maps $G_i$ is injective on the open set 
$*\times \mD^{n-2}_\circ$ (the spine).  Hence the
Strong Jordan separation tells us that each $G_i(S^{n-2})$ separates 
$S^{n-1}$, and that there are precisely two connected open sets
$U_i,V_i\subset S^{n-1}-G_i(S^{n-2})$  which
contain $G_i(\mD^{n-2})$ in their closure.  Furthermore, each of the
maps $F_i$ surjects onto either $V_i$ or $U_i$.

Now let $p:=\pi_x(\gamma(\infty))$,
and note that we have canonical isomorphisms:
$$H(p,G_1)\cong H(p,G_2)\cong H(p,G_3)\cong \mZ.$$
Indeed, this follows from the fact that all three maps have to coincide in 
a small enough neighborhood of $p$.  Now fix a generator $\{[x_m-y_m]\}$
for $H(p,G_1)$ as in Lemma 2.7.  Now Lemmas 2.5 and 2.6 imply that 
(upto re-indexing), we may assume that
for $m$ sufficiently large we have $x_m\in U_i$, $y_m\in V_i$ for all
$1\leq i\leq 3$.  We now assume, without loss of generality, that $F_1$ and $F_2$  
surject onto $V_1$ and $V_2$ respectively.  

Next we note that if $\epsilon$ is small enough, we will have that the ball $D$
of radius $\epsilon$ centered at $p$ only intersects $G_i(\mD^{n-2})$ (this 
follows from the argument in Claim 3).  In particular, each connected 
component of $D-G_i(\mD^{n-2})$ is contained in either $U_i$ or in $V_i$.
Now observe that the argument given in Lemma 2.5 ensures that there
are precisely two connected components of $D-G_i(\mD^{n-2})$ which contain
$p$ in their closure.  Denote these two components by $U$,$V$, labelled
so that for $m$ sufficiently large, $x_m\in U$, $y_m\in V$.  By the
choice of our labelling, we have that $U\subset \cap U_i$, and $V\subset
\cap V_i$.  

Finally, we note that $V$ is connected and contains $p$ in its closure
(by construction), and since $V\subset V_1\cap V_2$, we have that both
$F_1$ and $F_2$ surject onto $V$ (since they surject onto $V_1,V_2$ 
respectively).  This concludes the proof of Claim 5.

\vskip 10pt

\begin{Claim2}
Let $V\subset S^{n-1}$ be a connected open set, containing the 
point $p:=\pi_x(\gamma (\infty))$ in its closure.  Then $V$ contains
a connected open set $U_j$ lying in the complement of the set $I_x$.
\end{Claim2}

The proof of Claim 6 is identical to the one given in \cite{la}.  We reproduce
it here for completeness.
We first claim that the connected open set $V$ contains a point from 
$I_x$.  Indeed, if not, then $V$ would lie entirely in the complement
of $I_x$, hence would lie in some $U_i$.  Since $\pi_x(\gamma (\infty))$
lies in the closure of $V$, it would also lie in the closure of $U_i$,
contradicting the fact that $\gamma$ is {\it not} asymptotic to any of
the boundary components of the chamber containing $\gamma$.

So not only does $V$ contain the point $\pi_x(\gamma (\infty))$ in its
closure, it also contains some point $q$ in $I_x$.  We claim it in fact 
contains a point in $I_x-\cup (\partial U_i)$.  If $q$ itself lies in
$I_x-\cup (\partial U_i)$ then we are done.  The other possibility is 
that $q$ lies in the boundary of one of the $U_i$.  Now
since $V$ is connected, there exists a path $\eta$ joining $q$ to 
$\pi_x(\gamma (\infty))$.  Now assume that $\eta \cap \big(I_x-\cup 
(\partial U_i) \big)
=\{\pi_x(\gamma (\infty))\}$.  Let $\bar U_i$ denote the closed disks
(closure of the $U_i$),
and note that the complement of the set $I_x-\cup (\partial U_i)$ is the set $\cup (\bar U_i)$.

A result of Sierpinski \cite{sp} states the following: let
$X$ be an arbitrary topological space, $\{C_i\}$ a countable
collection of disjoint path connected closed subsets in $X$.  Then
the path connected components of $\cup C_i$ are precisely the
individual $C_i$. 

Applying Sierpinski's result to the set $\cup (\bar U_i)$ shows 
that the path connected component of this union are precisely the
individual $\bar U_i$.  So if $\eta \cap \big(I_x-\cup 
(\partial U_i) \big)=\{\pi_x(\gamma (\infty)\}$, we see that the
path $\eta$ must lie entirely within the $\bar U_i$ containing
$q$.  This again contradicts the fact that $\pi_x(\gamma (\infty))\notin
\cup (\partial U_i)$.
Finally, the fact that $V$ contains a point in $I_x-\cup (\partial U_i)$
allows us to invoke Claim 2, which tells us that there is some $U_j$ which 
is contained entirely within the set $V$, completing our argument for 
Claim 6.
\vskip 10pt

We now proceed to complete the proof of Proposition 3.2.  From Claims 5 and
6, there exists an open set $V\subset S^{n-1}=lk(x)$ with the property that:
\begin{itemize}
\item the maps $F_1$, $F_2$ surject onto $V$.
\item $V$ contains a connected open set $U_j$ lying in the complement of
the set $I_x$.  
\end{itemize}
This of course implies that both $F_1$ and $F_2$ must surject onto the 
boundary $\partial U_i$ of the $U_i$.  But now we observe that $\partial U_i$
lies in the set of injectivity $I_x$.  Since $\partial U_i$ lies in the
image of $F_1$ and $F_2$, the pre-image of $\partial U_i$ must lie in $L_1\cap L_2
=*\times \mD^{n-2}$.  Furthermore, the pre-image of $\partial U_i$ is homeomorphic
to an $S^{n-2}$.  So we obtain an embedded $S^{n-2}$ inside $*\times \mD^{n-2}$,
which immediately gives us a contradiction.  This completes the proof of Proposition 3.2.

\subsection{Connectedness and separation properties.}

We briefly state four results that will be used to deal with the second case
of Theorem 1.4, and will also be used in the arguments for Theorems
1.1 and 1.2.  The proofs of these results given in \cite{la} under the hypothesis
$n=3$ extend immediately to the more general setting $n\geq 2$.

\begin{Lem}
Let $\tilde B_i$ be a connected lift of the branching locus, and
let $\tilde W_j$, $\tilde W_k$ be two lifts of chambers which are
both incident to $\tilde B_i$.  Then $\tilde W_j- \tilde
B_i $ and $\tilde W_k- \tilde B_i $
lie in different connected components of $\tilde X^n - \tilde
B_i$.
\end{Lem}

\begin{Lem}
Let $\tilde W_i$ be a connected lift of a chamber, and let $\tilde
B_j$, $\tilde B_k$ be two connected lifts of the branching locus
which are both incident to $\tilde W_i$.  Then $\tilde B_j$ and
$\tilde B_k$ lie in different connected components of $\tilde X^n
- Int(\tilde W_i)$.
\end{Lem}

\begin{Lem}
Let $\partial^\infty \tilde B_i$ be the boundary at infinity of a
connected lift of the branching locus, and let $\partial^\infty
\tilde W_j$, $\partial^\infty \tilde W_k$ be the boundaries at
infinity of two lifts of chambers which are both incident to
$\tilde B_i$. Then $\partial^\infty \tilde W_j- \partial^\infty \tilde B_i $ and
$\partial^\infty \tilde W_k- \partial^\infty \tilde B_i$ lie in different 
connected components
of $\partial^\infty \tilde X^n - \partial^\infty \tilde B_i$.
\end{Lem}

\begin{Lem}
Let $\partial^\infty \tilde W_i$ be the boundary at infinity
corresponding to a connected lift of a chamber, and let
$\partial^\infty \tilde B_j$, $\partial^\infty \tilde B_k$ be the
boundary at infinity of two connected lifts of the branching locus
which are both incident to $\tilde W_i$. Then $\partial^\infty
\tilde B_j$ and $\partial^\infty \tilde B_k$ lie in different
connected components of $\partial^\infty \tilde X^n -
(\partial^\infty \tilde W_i - \cup \partial^\infty \tilde B_l)$, where the union 
is over all $\tilde B_l$ which are boundary components of $\tilde W_i$.
\end{Lem}

\subsection{The geodesic ray passes through infinitely many chambers.}

We now explain how to deal with the second case of the Theorem 1.4, namely we show the
following:

\begin{Prop}
Let $X^n$ be a simple, thick $n$-dimensional hyperbolic P-manifold, with $n\geq 3$.
Let $\gamma \subset \tilde X^n$ be a geodesic that passes through
infinitely many connected lifts $\tilde W_i$.  Then
$\gamma (\infty)$ is {\bf not} $(n-1)$-branching.
\end{Prop}

In the proof of Proposition 3.3, we will make use of a family of nice 
metrics on the boundary at infinity of an arbitrary $n$-dimensional hyperbolic 
P-manifold  (in fact, on the boundary at infinity of any 
CAT(-1) space).

\begin{Def}
Given an $n$-dimensional hyperbolic P-manifold, and a basepoint
$*$ in $\tilde X^n$, we can define a metric on the boundary at
infinity by setting $d_\infty (p,q)=e^{-d(*,\gamma_{pq})}$, where
$\gamma_{pq}$ is the unique geodesic joining the points $p,q$ (and $d$ denotes 
the distance inside $\tilde X_n$).
\end{Def}

The fact that $d_\infty$ is a metric
on the boundary at infinity of a proper $CAT(-1)$ space follows
from Bourdon (Section 2.5 in \cite{bo}).  Note that changing the basepoint from $*$
to $*^\prime$ changes the metric, but that for any $p,q\in
\partial^\infty(X^n)$, we have the inequalities:
$$A^{-1}\cdot d_{\infty, *}(p,q)\leq d_{\infty, *^\prime}(p,q)
\leq A\cdot d_{\infty, *}(p,q)$$ where $A=e^{d(*,*^\prime)}$, and
the subscripts on the $d_\infty$ refers to the choice of basepoint
used in defining the metric.   In particular, different choices
for the basepoint induce the same topology on $\partial ^\infty
\tilde X^n$, and this topology coincides with the standard topology
on $\partial ^\infty \tilde X^n$ (the quotient topology inherited from
the compact-open topology in the definition of $\partial ^\infty \tilde
X^n$ as equivalence classes of geodesic rays in $\tilde X^n$).
This gives us the freedom to select basepoints at our
convenience when looking for {\it topological} properties of the
boundary at infinity.  We now proceed to prove Proposition 3.3:

\begin{Prf}
The approach here consists of reducing to the situation covered
in proposition 3.2.
We start by re-indexing the various consecutive connected lifts
$\tilde W_i$ that $\gamma$ passes through by the integers.  Fix a
basepoint $x\in \tilde W_0$ interior to the connected lift $\tilde
W_0$, and lying on $\gamma$. Now assume that there is an injective
map $f:T\times \mD^{n-2}\longrightarrow
\partial ^\infty \tilde X^n$ with $\gamma (\infty)\in
f(*\times \mD^{n-2}_\circ)$. 

We start by noting that, between successive connected lifts
$\tilde W_i$ and $\tilde W_{i+1}$ that $\gamma$ passes through,
lies a connected lift of the branching locus, which we denote
$\tilde B_i$.  Observe that distinct connected lifts of the
branching locus stay a uniformly bounded distance apart. Indeed,
any minimal geodesic joining two distinct lifts of the branching
locus must descend to a minimal geodesic in a $W_i$ with endpoints
in the branching locus.  But the length of any such geodesic is
bounded below by the injectivity radius of $DW_i$, the double of 
$W_I$ across its boundary.  By setting
$\delta$ to be the infimum, over all the finitely many chambers
$W_i$, of the injectivity radius of the doubles $DW_i$, we have $\delta>0$.
Let $K_i$ be the connected component of $\partial ^\infty \tilde
X^n -\partial ^\infty \tilde B_i$ containing $\gamma (\infty)$.
Then for every $p\in K_i$ ($i\geq 1$), we have:
$$d_x(p, \gamma (\infty))<e^{-\delta (i-1)}.$$
Indeed, by Lemma 3.1, $\tilde
B_i$ separates $\tilde X^n$ into (at least) two totally geodesic
components. Furthermore, the component containing $\gamma(\infty)$
is {\it distinct} from that containing $x$.  Hence, the distance
from $x$ to the geodesic joining $p$ to $\gamma (\infty)$ is at 
least as large as the
distance from $x$ to $\tilde B_i$.  But the later is bounded below
by $\delta (i-1)$.  Using the definition of the metric at infinity, and picking 
$x$ as our basepoint, our estimate follows.  

Since our estimate shrinks to zero, and since the distance from
$\gamma (\infty)$ to $f(\partial (T\times \mD^{n-2}))$ is positive, we 
must have a
point $q\in f(*\times \mD^{n-2}_\circ)$ satisfying $d_x(q,
\gamma(\infty))>e^{-\delta (i-1)}$
for $i$ sufficiently large. Since $\partial ^\infty \tilde B_i$ 
separates (by Lemma 3.3), we see that for $i$ sufficiently
large, $f(*\times \mD^{n-2}_\circ)$ contains points in two distinct components
of $\partial ^\infty \tilde X^n-\partial ^\infty \tilde B_i$.  
This implies that there
is a point $q^\prime \in f(*\times \mD^{n-2}_\circ)$ that lies within 
some $\partial ^\infty
\tilde W_k-\cup\partial^\infty \tilde B_l$, where the union runs over all
$\tilde B_l$ which are boundary components of $\tilde W_k$. 
But such a point corresponds to a geodesic ray lying
entirely within $\tilde W_k$, and {\it not} asymptotic to any of
the lifts of the branching locus.  Finally, we note that {\it any}
point in the image $f(*\times \mD^{n-2}_\circ)$ can be considered
$(n-1)$-branching, so in particular the point $q^\prime$ is
$(n-1)$-branching.  But in Proposition 3.2, we showed
this is impossible.  Our claim follows.
\end{Prf}

Finally, combining Proposition 3.1, 3.2, and 3.3 gives us Theorem 1.4.

\section{Rigidity results.}

In this section, we discuss two types of rigidity results (Theorems 1.1 and 1.2
from the introduction) for simple, thick
hyperbolic P-manifolds of dimension $\geq 3$.  In \cite{la}, the author proved
both these theorems for the case $n=3$, by making use of the $3$-dimensional
analogue of Theorem 1.4.  But the argument given in that previous paper works
equally well in higher dimensions (once one knows that Theorem 1.4 holds).
As such, we merely outline
the arguments in this section, and refer the interested reader to \cite{la} 
for more details.  We start by
recalling a useful result from \cite{la}:

\begin{Lem}
Let $X^n$ be a simple hyperbolic P-manifold of dimension
at least three
Let $\partial ^\infty \tilde B\subset \partial ^\infty \tilde X^n$
consist of all limit points of geodesics in the branching locus.
If $n\geq 3$, then the maximal path-connected components of
$\partial ^\infty \tilde B$ are precisely the sets of the form
$\partial ^\infty \tilde B_i$, where $\tilde B_i\subset \tilde B$ is a
single connected component of the lifts of the branching locus.
\end{Lem}

Indeed, each $\partial ^\infty \tilde B_i$ is path-connected (since
$n\geq 3$ and each $\tilde B_i$ is isometric to $\mH^{n-1}$) 
and closed (since each $\tilde B_i$ is totally geodesic), hence
$\partial ^\infty \tilde B$ is the union of a disjoint family of closed,
path-connected subspaces.  The fact that the maximal path-connected 
components of $\partial ^\infty \tilde B$ are the individual 
$\partial ^\infty \tilde B_i$ now follows from an easy
application of Sierpinski's result \cite{sp}.

Note that in the previous Lemma, the hypothesis that the dimension be 
$\geq 3$ really is necessary, since in the $2$-dimensional case one would
have the codimension one strata consisting of geodesics.  Hence for 
each connected component of a lift of the branching locus, we would have
that $\partial ^\infty \tilde B_i$ consists of a pair of points (and
hence, would not be path connected).

\subsection{Mostow rigidity.}

We now explain how Theorem 1.1 follows from Theorem 1.4.  Theorem 1.1 is 
clearly an analogue of Mostow rigidity for simple, thick hyperbolic P-manifolds.  
Let us start with a pair $X_1, X_2$ of simple, thick hyperbolic P-manifolds 
of dimension 
$\geq 3$, and an isomorphism $\phi:\pi_1(X_1)\rightarrow \pi_1(X_2)$ on the
level of the fundamental groups.  Our goal is to show that this implies that
$X_1$ is isometric to $X_2$.  

In order to do this, we note that an abstract isomorphism between the 
fundamental groups induces a quasi-isometry between the universal covers 
$\tilde \phi:\tilde X_1\rightarrow \tilde X_2$.  
One of the basic facts about quasi-isometries between $\delta$-hyperbolic
spaces is that they extend to give {\it homeomorphisms} between the boundaries
at infinity of the spaces.  In particular, the quasi-isometry $\tilde \phi
:\tilde X_1\rightarrow \tilde X_2$ extends to give a homeomorphism, denoted
$\phi_\infty$ from $\partial ^\infty \tilde X_1$ to $\partial ^\infty \tilde
X_2$.  Furthermore, this homeomorphism intertwines the actions of the
groups $\pi_1(X_1)$ and $\pi_1(X_2)$ on the respective boundaries at infinity.

Now recall that the property of being 
$(n-1)$-branching is a topological invariant, $\phi_\infty$ must
map branching points in $\partial ^\infty \tilde X_1$ to branching points
in $\partial ^\infty \tilde X_2$.  Since connectedness is preserved under
homeomorphisms, $\phi_\infty$ must
map connected components of the set of branching points in 
$\partial ^\infty \tilde X_1$ to connected components of the set of
branching points in $\partial ^\infty \tilde X_2$.  

Now Theorem 1.4 and Lemma 4.1 combine to tell us precisely what these connected 
components are:
each connected component corresponds precisely to the points in the boundary at 
infinity of a single connected lift of the branching locus. 
This forces $\phi_\infty$ to induce a bijection between the lifts of the 
branching locus in $\tilde X_1$ and the lifts of the branching locus in 
$\tilde X_2$.

The separation argument (Lemmas 3.3 and 3.4) in the previous section ensure 
that the map 
$\phi_\infty$ takes the boundaries of the lifts of chambers in $\partial 
^\infty \tilde X_1$ to the boundaries of the lifts of chambers in $\partial 
^\infty \tilde X_2$.  This naturally induces a bijection between the lifts
of chambers in $\tilde X_1$ and the lifts of chambers in $\tilde X_2$.  
Furthermore, we note that this bijection preserves the adjacency relation
between chambers (again, this follows from Lemmas 3.3 and 3.4).

Finally, we note that the fundamental group of the chambers in a simple,
thick hyperbolic P-manifold $X$ can be detected
from the boundary at infinity: indeed, the fundamental group of a chamber $W_i$
is isomorphic to the stabilizer of $\partial ^\infty \tilde W_i$ under the
action of $\pi_1(X)$ on $\partial ^\infty \tilde X$.  Since the map $\phi_\infty$
intertwines the actions, we conclude that the bijection takes lifts of
chambers to lifts of chambers with isomorphic fundamental groups.  

Mostow rigidity for hyperbolic manifolds with totally geodesic boundary (see
Frigerio \cite{fr}) now allows us to conclude that the individual chambers must
be isometric.  Furthermore, the gluings between the chambers are detected, since lifts 
of adjacent chambers
map to lifts of adjacent chambers, and the dynamics of $\pi_1(X)$ on $\partial 
^\infty \tilde X$ allow us to identify the subgroups of the individual chambers 
(i.e. the boundary components) that get identified together.  Patching this
information together, we obtain an equivariant isometry from $\tilde X_1$ to
$\tilde X_2$.  This naturally descends to an isometry from $X_1$ to $X_2$.  Furthermore,
by construction this isometry must induce $\phi$ on the level of $\pi_1$, concluding
the argument for Theorem 1.1.

\subsection{Quasi-Isometric rigidity.}

In this section, we outline the argument for showing the quasi-isometric rigidity
result (Theorem 1.2).  The starting point is the following well known result
(for a proof, see for instance Prop. 3.1 in \cite{fb}):

\begin{Lem}
Let $X$ be a proper geodesic metric space, and assume that every
quasi-isometry from $X$ to itself is in fact a bounded distance (in the sup norm)
from an isometry. Furthermore, assume that a finitely generated
group $G$ is quasi-isometric to $X$.  Then there exists a
cocompact lattice $\Gamma\subset Isom(X)$, and a finite group $F$
which fit into a short exact sequence:
$$0\longrightarrow F\longrightarrow G\longrightarrow \Gamma
\longrightarrow 0$$
\end{Lem}

In particular, if one can show that on the universal cover $\tilde X^n$ of a 
simple, thick hyperbolic P-manifold of dimension $n\geq 3$, every quasi-isometry
is a bounded distance from an isometry, then applying Lemma 4.2 would yield 
Theorem 1.2.

So let us assume that $\phi:\tilde X^n\rightarrow \tilde X^n$ is a quasi-isometry,
and observe that this quasi-isometry induces a homeomorphism $\partial^\infty \phi:
\partial ^\infty \tilde X^n\rightarrow \partial ^\infty \tilde X^n$.  Arguing as
in the previous section, this quasi-isometry must map each $\partial ^\infty \tilde W_i$
($W_i$ a chamber in $X^n$) homeomorphically to $\partial ^\infty \tilde W_j$
(where $W_j$ is again a chamber in $X^n$).  In particular, the quasi-isometry $\phi$
restricts (upto a bounded amount) to a quasi-isometry from $\tilde W_i$ to $\tilde W_j$.
We now recall the well-known:

\vskip 5pt

\noindent {\bf Folklore Theorem:}
{\it Let $M^n_1$, $M^n_2$ be a pair of $n$-dimensional ($n\geq 3$) hyperbolic manifolds 
with non-empty, totally geodesic boundaries.  If $\phi:\tilde M^n_1\rightarrow
\tilde M^n_2$ is a quasi-isometry, then there exists an isometry $\psi:\tilde M^n_1
\rightarrow \tilde M^n_2$ at a bounded distance from $\phi$ (in the sup norm).}

\vskip 5pt

Proofs of this result have been discovered at various times by
Farb, Kapovich, Kleiner, Leeb, Schwarz,
Wilkinson, and others, though to the author's knowledge the only written proof 
of this theorem is to be found in Frigerio \cite{fr2}.

By appealing to the previous
Theorem, we see that the quasi-isometry of $\tilde X^n$ has the property that, for 
each lift of a chamber of $\tilde X^n$, there is an isometric mapping to 
the lift of some other chamber that
is bounded distance away from the original quasi-isometry.  There are two points to
verify: 
\begin{enumerate}
\item that the individual isometries on the lifts of the chambers glue 
together to give a global isometry,
\item that the distance between the quasi-isometry and the isometry on each lift of
a chamber is uniformly bounded.
\end{enumerate}
Concerning point (1), we first note that the quasi-isometry maps adjacent lifts of 
chambers
to adjacent lifts of chambers (this follows from Lemmas 3.3 and 3.4).  The fact that the
isometries coincide on the common boundary follows from the fact that the 
isometries induced by the incident chambers are all at a bounded distance from 
each other (since they are all at a bounded distance from the original quasi-isometry).  
But the common boundary is isometric to $\mH^{n-1}$ (by the simplicity hypothesis), 
and any two
isometries of $\mH^{n-1}$ which are bounded distance apart have to coincide.  This
allows us to glue together the isometries, taking care of point (1).

Concerning point (2), one sees that it is possible, on the universal covers of  
compact hyperbolic
manifolds with totally geodesic boundary, to give a uniform bound (depending only
on the quasi-isometry constants) on the distance between a quasi-isometry and 
the isometry that is at finite distance from it (see the end of Section 3.2 in
\cite{la}).  Since the quasi-isometries on
all the chambers were induced by a global quasi-isometry of $\tilde X^n$, they
all have the same quasi-isometry constants, and hence one obtains point (2).
This completes our discussion of Theorem 1.2.

\section{Concluding remarks.}

It would be interesting to find some further applications of the rigidity results
presented in this paper.  One potential application is to the Gromov-Thurston 
examples of negatively curved manifolds.  Recall that these manifolds $\bar M^n$ are 
topologically
ramified coverings of hyperbolic manifolds $M^n$, with the ramification occuring over a 
totally geodesic codimension two submanifold $N^{n-2}$ (see \cite{gt}).  These manifolds
naturally inherit a (singular) CAT(-1) metric from the hyperbolic metric on $M^n$, with
the singular set for the metric consisting of the pre-image of the ramification
locus $N^{n-2}$.  Gromov-Thurston showed that this singular CAT(-1) metric can be smoothed
out to a negatively curved Riemannian metric.

Now consider the case where $N^{n-2}$ bounds a totally geodesic codimension one 
submanifold $K^{n-1}$, and observe that the pre-image of $K^{n-1}$ 
form a totally geodesic subset in $\bar M^{n}$ (equipped with the singular CAT(-1)
metric) isometric to a simple, thick (assuming the ramified cover has degree at 
least three) hyperbolic P-manifold.  It is tempting to see whether the rigidity 
results presented in this paper can be used to analyze the Gromov-Thurston examples.
One basic questions along these lines is the following:

\vskip 5pt

\noindent {\bf Question:} Let $M^n$ be a Gromov-Thurston negatively curved
manifold, $N^{n-2}\subset M^n$ the ramification locus, 
and $\Gamma:=\pi_1(M^n)$, $\Lambda :=\pi_1(N^{n-2})$ their respective fundamental
groups.  If $\phi:\Gamma\rightarrow \Gamma$
is an abstract isomorphism, is it true that $\phi(\Lambda)$ is always 
conjugate to $\Lambda$?  

\vskip 5pt

Next, we note that the argument for Theorem 1.4 (Characterization of branching 
points) extends (with suitable modifications in the proofs) to the more general
setting of P-manifolds with {\it negatively curved} metrics.  We did not pursue
this in the present paper because we could not think of suitable applications 
of this result.  Of course, there is no hope of proving Mostow type rigidity
in this setting (though topological rigidity \`a la Farrell-Jones \cite{fj}
should still hold).  As for quasi-isometric rigidity, what would be needed is
the following open question:

\vskip 5pt

\noindent {\bf Question:} If $M^n$ is a negatively curved Riemannian manifold 
with non-empty, totally geodesic boundary.  Is every quasi-isometry of the 
universal cover $\tilde M^n$ bounded distance away from an isometry?

\vskip 5pt

The argument for the Folklore Theorem (see \cite{fr2}) 
cannot possibly work in this more general setting, since it relies heavily on the
relationship between quasi-isometries (and isometries) of $\mH^n$ on the one hand, and
quasi-conformal (and conformal) maps of $S^{n-1}=\partial ^\infty \mH^n$ on the other.
Since no such correspondance exists in the general variable curvature case
(though see Pansu \cite{pa}), an
entirely new approach would need to be given.  Finally, concerning Theorem 1.1, we
ask the following:

\vskip 5pt

\noindent {\bf Question:} If one has a map $f:S^n\rightarrow S^{n+1}$ with the
property that the set of injectivity has {\it positive measure}, does that imply
that $f(S^n)$ separates $S^{n+1}$?

\vskip 5pt

It is clear that the approach we use really requires the presence of an open disk in
the set of injectivity, and hence is of no use in answering this last question.

\end{document}